\documentclass{gtart}
\usepackage{amsmath, amssymb, stmaryrd, verbatim}
\usepackage{tipa}
\usepackage{graphs}
\usepackage{graphicx}
\newtheorem{thm}{Theorem}[section]
\newtheorem{cor}[thm]{Corollary}
\newtheorem{lem}[thm]{Lemma}

\theoremstyle{definition}
\newtheorem{defn}[thm]{Definition}

\theoremstyle{remark}

\numberwithin{equation}{section}

\newcommand{\al}{\alpha}
\newcommand{\be}{\beta}
\newcommand{\de}{\delta}

\newcommand{\ga}{\gamma}

\newcommand{\Ga}{\Gamma}

\newcommand{\Om}{\Omega}
\newcommand{\om}{\omega}

\newcommand{\Z}{\mathbb Z}
\newcommand{\Q}{\mathbb Q}

\newcommand{\del}{\partial}

\newcommand{\arc}{\texttoptiebar}

\DeclareMathOperator{\fr}{\mathrm{frac}}

\begin{document}
\title{Stein fillable Seifert fibered 3--manifolds}

\author{Ana G.~Lecuona, Paolo Lisca}
\address{Dipartimento di Matematica ``L. Tonelli'',\\
Universit\`a di Pisa \\
Largo Bruno Pontecorvo 5\\
I-56127 Pisa, Italy} 
\email{lecuona@mail.dm.unipi.it}
\secondemail{lisca@dm.unipi.it}


\begin{abstract}
We characterize the closed, oriented, Seifert fibered 3--manifolds 
which are oriented boundaries of Stein manifolds. We also show that 
for this class of 3--manifolds the existence of Stein fillings is equivalent to 
the existence of symplectic fillings. 
\end{abstract}
\maketitle

\section{Introduction and statement of results}\label{s:intro}

The most important dichotomy in 3--dimensional contact topology is the one introduced 
by Eliashberg between tight and overtwisted contact structures (see e.g.~\cite{Et, Ge}). 
Nowadays there are several different ways to prove that a contact structure is tight, 
but for a long time the only systematic way to construct tight contact structures on a 
closed 3--manifold $Y$ was to show $Y$ to be orientation preserving diffeomorphic 
to the oriented boundary of a Stein manifold and then appeal to a theorem of Eliashberg 
and Gromov~\cite{El2, Gr}. This naturally led to the question of which 3--manifolds 
carry tight contact structures, as well as to the related question of which 3--manifolds 
admit \emph{Stein fillings}, i.e.~are orientation preserving 
diffeomorphic to the boundary of a Stein manifold. 
The first example of an oriented 3--manifold admitting no Stein fillings was provided in~\cite{Li98}, 
and infinitely many examples were found in~\cite[Theorem~4.2]{LS04-2}  and~\cite[Proposition~4.1]{LS07-1}. 
While the classification of the closed, Seifert fibered 3--manifolds carrying 
tight contact structures was recently achieved~\cite{LS09}, the classification of the  
Stein fillable ones was still missing. The purpose of the 
present paper is to fill this gap. Our main result, Theorem~\ref{t:main} below, 
identifies explicitely the family of closed, oriented, Seifert fibered 3--manifolds 
which are orientation preserving diffeomorphic to the boundary of a Stein manifold.

\begin{figure}[ht]
\centering
\includegraphics[scale=0.5]{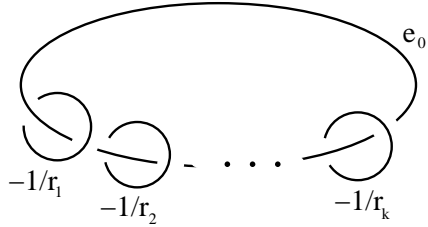}
\caption{The Seifert 3--manifold $Y(e_0;r_1,\ldots,r_k)$}
\label{f:seifert}
\end{figure}

We need some preliminaries in order to state our results.  
Eliashberg~\cite{El1} proved that smooth, even--dimensional manifolds 
carrying Stein structures can be characterized as having suitable handle 
decompositions. Gompf~\cite[Theorem~5.4]{Go} elaborated on 
Eliashberg's result to show that a closed, oriented, Seifert fibered 
3--manifold $Y$ admits a Stein filling unless $Y$ is orientation preserving 
diffeomorphic to the oriented 3--manifold $Y(e_0;r_1,\ldots,r_k)$ given by 
the surgery description of Figure~\ref{f:seifert}, where 
\[
e_0=-1,\quad k\geq 3\quad\text{and}\quad 1>r_1\geq \cdots \geq r_k>0.
\] 
Gompf also discovered a sufficient condition for the existence of Stein fillings 
of a 3--manifold  of the form $Y(-1;r_1,\ldots,r_k)$. We describe and use this condition 
in Section~\ref{s:existence}. 

\begin{defn}\label{d:realizable} 
A $k$--tuple $(r_1,\ldots,r_k)\in (\Q\cap (0,1))^k$ with $k\geq 3$ and 
$r_1\geq r_2\geq \cdots \geq r_k$ is \emph{realizable} if 
there exist coprime integers $n>h>0$ such that 
\[
\frac{h}{n} > r_1,\quad \frac{n-h}n > r_2,\quad\text{and}\quad \frac1n > r_3,\ldots,r_k.
\]
\end{defn}

\begin{defn}\label{d:specialtype}
A closed, oriented, Seifert fibered 3--manifold is~\emph{of special type} 
if it is orientation preserving diffeomorphic to $Y(-1;r_1,\ldots, r_k)$, where $k\geq 3$, 
$1>r_1\geq r_2\geq \cdots \geq r_k>0$ and the following conditions both hold: 
\begin{enumerate}
\item
$(r_1,\ldots,r_k)$ is not realizable;
\item
$r_1+\cdots + r_k > 1 > r_1+r_2$.
\end{enumerate}
\end{defn}

In Section~\ref{s:existence} we use Gompf's sufficient condition 
for the existence of Stein fillings of $Y(-1;r_1,\ldots, r_k)$ to establish the following: 

\begin{thm}\label{t:mainexistence} 
Let $Y$ be a closed, oriented, Seifert fibered 3--manifold which is not of special type. 
Then, $Y$ is orientation preserving diffeomorphic to the boundary of a Stein surface. 
\end{thm} 

Recall (see e.g.~\cite{Ge}) that a \emph{symplectic filling} of a contact 3--manifold $(Y,\xi)$ 
is a pair $(X,\om)$, where $X$ is a smooth 4--manifold with boundary oriented by 
a symplectic form $\om\in\Om^2(X)$, and such that there is an orientation preserving 
diffeomorphism $\varphi\co Y\to\del X$ with $\om|\varphi_*(\xi)\neq 0$ at each 
point of $\del X$. A Stein filling is a symplectic filling but the converse is not true, 
because there are examples of symplectically fillable contact 3--manifolds 
which are not Stein fillable~\cite{Gh}. Similarly, there exist several examples of 
tight, contact Seifert fibered 3--manifolds which are not symplectically 
fillable~\cite{LS03, LS04-1,LS04-2,GLS}. 
In Section~\ref{s:nonexistence} we apply Donaldson's theorem on 
the intersection forms of definite 4--manifolds to prove the following: 

\begin{thm}\label{t:mainnonexistence}
A closed, oriented, Seifert fibered 3--manifold of special type admits 
no symplectic fillings. 
\end{thm}

Combining Theorems~\ref{t:mainexistence} and~\ref{t:mainnonexistence}  
immediately gives our main result:

\begin{thm}\label{t:main}
Let $Y$ be a closed, oriented, Seifert fibered 3--manifold. Then, the 
following conditions are equivalent:
\begin{enumerate}
 \item 
$Y$ admits Stein fillings;
\item 
$Y$ admits symplectic fillings;
\item
$Y$ is not of special type.
\end{enumerate}
\end{thm}

\begin{proof}
A Stein filling is a symplectic filling, therefore (1) implies (2). 
By Theorem~\ref{t:mainnonexistence} (2) implies (3) and by 
Theorem~\ref{t:mainexistence} (3) implies (1). 
\end{proof}

Theorem~\ref{t:main} implies Corollary~\ref{c:main} below, which shows that 
the condition that a 3--manifold is Seifert fibered of special type can be 
reformulated in terms of open book decompositions. Recall that an 
open book decomposition of a 
closed 3--manifold is called~\emph{positive} if its monodromy
can be written as a product of right--handed Dehn twists.
Loi and Piergallini~\cite[Theorem~4]{LP01} proved that a smooth, closed, 
oriented 3--manifold $Y$ is the boundary of a Stein surface if and only if $Y$ 
admits a positive open book decomposition. Combining this result with 
Theorem~\ref{t:main} immediately yields the following:

\begin{cor}\label{c:main}
A closed, oriented, Seifert fibered 3--manifold $Y$ admits a positive open book 
decomposition if and only if $Y$ is not of special type.  
\qed\end{cor}

We would like to end this introduction with a few remarks about ``how special'' 
Seifert 3--manifolds of special type are. First, it should be clear 
from the definition that there exist infinitely many oriented Seifert 3--manifolds 
of special type. Indeed, infinitely many examples of closed, oriented 
Seifert fibered 3--manifolds without symplectic fillings are 
known (\cite[Theorem~4.2]{LS04-2}  and~\cite[Proposition~4.1]{LS07-1}).
According to Theorem~\ref{t:main} such examples must be all of special type, 
and in fact it can be easily verified that they are of special type. The infinitely many oriented 
Seifert 3--manifolds which do not carry tight contact structures~\cite{LS07-1, LS09} are also 
of special type. Second, if $Y=Y(-1;r_1,\ldots, r_k)$ is of special type then $-Y$ is not, because 
it is of the form $Y(-1+k;1-r_k,\ldots,1-r_1)$. This is consistent with the 
general fact that each oriented Seifert fibered 3--manifold admits a Stein filling after possibly 
reversing its orientation~\cite[Corollary~5.5(a)]{Go}. 
Also, Condition (1) from Definition~\ref{d:specialtype} 
together with the results of~\cite{LS07-2} imply that an oriented Seifert fibered 3--manifold of 
special type is an $L$--space in the sense of~\cite[Definition~1.1]{OS05}. 
(See~\cite{LS07-2} for several different characterizations of 
Seifert fibered $L$--spaces). In fact, it follows from~\cite{LS07-2} that if $Y$ is an oriented 
Seifert fibered 3--manifold which is an $L$--space, then after possibly reversing  
orientation $Y$ is of the form $Y(e_0;r_1,\ldots,r_k)$, with $k\geq 3$, 
$1>r_1\geq\cdots\geq r_k>0$ and either (i) $e_0\geq 0$ or 
(ii) $e_0=-1$ and $(r_1,\ldots,r_k)$ is not realizable. Therefore, the oriented 
Seifert fibered 3--manifolds of special type are precisely the oriented, Seifert 
fibered $L$--spaces of the form $Y(-1;r_1,\ldots,r_k)$ with $k\geq 3$, 
$1>r_1\geq\cdots\geq r_k>0$ and $(r_1,\ldots,r_k)$ satisfying Condition (2) 
from Definition~\ref{d:specialtype}.  

The organization of the paper is straightforward: in Section~\ref{s:existence} 
we prove Theorem~\ref{t:mainexistence} and in 
Section~\ref{s:nonexistence} we prove Theorem~\ref{t:mainnonexistence}. 

\section{Existence of Stein fillings}\label{s:existence}

The purpose of this section is to prove Theorem~\ref{t:mainexistence}. We start with  
recalling Gompf's sufficient condition from~\cite{Go} for the existence of a Stein filling 
of $Y(-1;r_1,\ldots,r_k)$. 

Given a rational number $r\in\Q$ we define an integer $\llbracket r\rrbracket\in\Z$ 
by setting $r=\llbracket r\rrbracket + \fr(r)$, where $\fr(r)\in [0,1)$. 
Define 
\[
r'_i:=-\frac1{r_i},\quad i=1,\ldots,k. 
\]
Let $s\in (-\infty,-1)$ be such that $\frac1s:=-1-\frac{1}{r'_1}$. If $s\neq r'_2$ then 
it is easy to check that there is a map 
\[
 A\co\Q\cup\{\infty\}\to\Q\cup\{\infty\}
\]
of the form $A(r)=\frac{c+dr}{a+br}$, such that: 
\begin{equation}\label{e:g1}
ad-bc=\pm 1,
\end{equation}
\begin{equation}\label{e:g2}
A(s)\in (-1,0],
\end{equation}
\begin{equation}\label{e:g3}
A(r'_2)\in [-\infty,-1).
\end{equation}
Let  
\[
 t:=
\begin{cases} 
 0\quad\text{if}\ A(0)\in [0,+\infty]\\
\frac1{A(s)}\quad\text{if}\ A(0)\in [-1,0)\\
A(r'_2)\quad\text{if}\ A(0)\in (-\infty,-1).
\end{cases}
\]
Set $M:=\max(|a|,|c|)$, $m:=\min(|a|,|c|)$ and 
\[
 n_A(r'_1,r'_2):=-m(\llbracket t\rrbracket + 1) - M.
\]
Finally, let 
\[
n(r'_1,r'_2):=
\begin{cases}
0 & \quad\text{if}\ s = r'_2,\\
\sup_A n_A(r'_1,r'_2) & \quad\text{if}\ s\neq r'_2.
\end{cases}
\]
Gompf~\cite{Go} shows that $Y(-1;r_1,\ldots,r_k)$ 
is the boundary of a Stein surface if:
\begin{equation}\label{e:Gcond}
 n(r'_1,r'_2)>r'_3,\ldots,r'_k.
\end{equation}
Observe that when $s=r'_2$ Condition~\eqref{e:Gcond} is 
automatically satisfied.

In order to prove Theorem~\ref{t:mainexistence} we need two results. The first 
result is Theorem~\ref{t:realizable}, which establishes the existence of a 
Stein filling for $Y(-1;r_1,\ldots,r_k)$ 
under the assumption that the $k$--tuple $(r_1,\ldots,r_k)$ is~\emph{realizable}, 
that is to say that there exist coprime integers $n>h>0$ such that 
$r_1<\frac{h}{n}$, $r_2<\frac{n-h}n$ and $r_3,\ldots,r_k<\frac1n$. 

\begin{thm}\label{t:realizable}
Suppose that $k\geq 3$, $1>r_1\geq r_2\geq \cdots \geq r_k>0$ and 
$(r_1,\ldots,r_k)$ is realizable. Then, $Y(-1;r_1,\ldots,r_k)$ is 
orientation preserving diffeomorphic to the boundary of a 
Stein surface. 
\end{thm}

\begin{proof} 
Recall that we defined $r'_i:=-\frac1{r_i}$, $i=1,\ldots, k$. 
We will prove that there is a map $A\co\Q\cup\{\infty\}\to\Q\cup\{\infty\}$
satisfying Properties~\eqref{e:g1}, \eqref{e:g2}, \eqref{e:g3} above, and such that 
\[
n_A(r'_1,r'_2)>r'_3,\ldots r'_k.
\] 
In view of Gompf's condition~\eqref{e:Gcond}, this clearly suffices to prove
the statement. 

By the realizability assumption, there is a positive 
integer $n_0$ such that, for some integer 
$h_0$ coprime with $n_0$ and satisfying 
$n_0>h_0>0$ we have
\[
r'_2<-\frac{n_0}{h_0}<s\quad\text{and}\quad 
-n_0 >r'_3,\ldots r'_k.
\]
Denote by $n$ the smallest positive integer such that, for 
some integer $h$ coprime with $n$ with $n>h>0$, we have 
\[
n\leq n_0\quad\text{and}\quad
r'_2<-\frac{n}{h}<s.
\]
Notice that, since $n\leq n_0$, $-n\geq -n_0>r'_3,\ldots,r'_k$. 
Moreover, $n$ and $h$ being coprime, there exist $a, b\in\Z$ such that 
\begin{equation}\label{e:coprime}
 a h - b n = 1.
\end{equation} 
If the pair $(a,b)$ solves Equation~\eqref{e:coprime}, so does the pair
$(a+zn,b+zh)$ 
for each $z\in\Z$. Therefore, we can choose a solution $(a,b)$ such that 
$0\leq a < n$. Indeed, since $a=0$ would imply $n=1$, which is not the 
case because $n>h>0$, we can assume $0<a<n$. From Equation~\eqref{e:coprime} 
we get
\[
 b = \frac{ah}{n} - \frac{1}{n},
\]
hence $-1/n<b<h-1/n$, which is equivalent to 
\[
 0\leq b < h.
\]
Soon it will be convenient to have $b>0$, therefore we deal now with 
the special case $b=0$. By Equation~\eqref{e:coprime}, $b=0$  
implies $a=h=1$, therefore $r'_2<-n<s$. Moreover, by the minimality of $n$ we
must 
have $-(n-1)\geq s$. Define the map $A$ by 
\[
 A(r):= r + n-1.
\]
This map is of the form $\frac{c+dr}{a+br}$ with $a=1$, $b=0$, $c=n-1$ and
$d=1$, 
therefore $ad-bc=1$. Clearly $A$ is monotone increasing, $A(-n)=-1$ and
$A(-n+1)=0$. Therefore 
$A(r'_2)\in (-\infty,-1)$ and $A(s)\in (-1,0]$. Therefore $A$ satisfies the
required Properties~\eqref{e:g1}, \eqref{e:g2} and \eqref{e:g3}. Since $A(0)=n-1\in [0,+\infty]$ 
we have $t=0$, $m=1$ and $M=n-1$, therefore 
\[
 n_A(r'_1,r'_2) = - m - M = -1 -n + 1 = -n>r'_3,\ldots,r'_k.
\]
From now on we assume $0<b<h$. Observe that Equation~\eqref{e:coprime} is 
equivalent to  
\[
 \frac{a}{b} = \frac{n}{h} + \frac{1}{hb}, 
\]
which implies 
\[
 -\frac{a}{b} < -\frac{n}{h}.
\]
In fact, by our choice of $n$ we must have 
\[
 -\frac{a}{b} \leq r'_2.
\]
Now we define $A$ by 
\[
 A(r):=\frac{(n-a)+(h-b)r}{a+br} = -1 + \frac{n+hr}{a+br}.
\]
For this map we have $c=n-a$ and $d=h-b$, therefore 
\[
ad-bc=a(h-b)-b(n-a)=ah-bn=1.
\] 
The map $A$ is monotone increasing for every $r\neq -\frac{a}{b}$, because 
\[
\frac{dA}{dr}(r) = \frac{1}{(a+br)^2}.
\]
Equation~\eqref{e:coprime} implies $-\frac{n-a}{h-b}>-\frac{n}{h}$, thus 
by the choice of $n$ we have $s\leq -\frac{n-a}{h-b}$. We conclude  
\[
 A(-\frac{a}{b})=-\infty \leq A(r'_2) < A(-\frac{n}{h})=-1 < A(s) \leq
A(-\frac{n-a}{h-b})=0.
\]
Therefore $A$ satisfies Properties~\eqref{e:g1}, \eqref{e:g2} and \eqref{e:g3}. Since 
\[
 A(0)= -1 + \frac{n}{a} = \frac{n-a}{a} \in (0,+\infty)
\]
we have $t=0$, thus 
\[
n_A(r'_1,r'_2)=-m-M=-|a|-|c|=-a-(n-a) = -n>r'_3,\ldots,r'_k.
\]
\end{proof} 

We can now move on to the second result needed for the proof of Theorem~\ref{t:mainexistence}, 
that is Theorem~\ref{t:r1r2big} below. This result will establish the existence of 
a Stein filling for $Y(-1;r_1,\ldots,r_k)$ under the assumption $r_1+r_2\geq 1$. 
The proof of Theorem~\ref{t:mainexistence} will then follow combining 
Theorems~\ref{t:realizable} and~\ref{t:r1r2big}.  

\begin{figure}[!ht]
\centering
\includegraphics[scale=0.35]{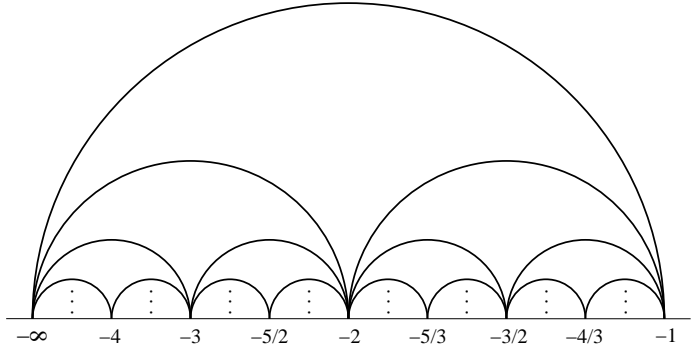}
\caption{Arcs of the Farey tessellation}
\label{f:farey}
\end{figure}

Consider the standard Farey tessellation of the hyperbolic plane. 
Figure~\ref{f:farey} illustrates some of the arcs of the  
tessellation with both endpoints in the interval $[-\infty,-1]$. 
We shall refer to any such arc with endpoints $\al<\be$ as to 
the~\emph{Farey arc} \arc{$\al\be$}. 

Observe that, given a Farey arc \arc{$\al\ga$}, there is a unique point 
$\be$ such that $\al<\be<\ga$ and there exist Farey arcs \arc{$\al\be$} and \arc{$\be\ga$}. 
In what follows, we shall refer to the unique point $\be$ as to the 
\emph{middle point} of \arc{$\al\ga$}, and denote it by $m(\al,\ga)$. 

\begin{figure}[!ht]
\centering
\includegraphics[scale=0.5]{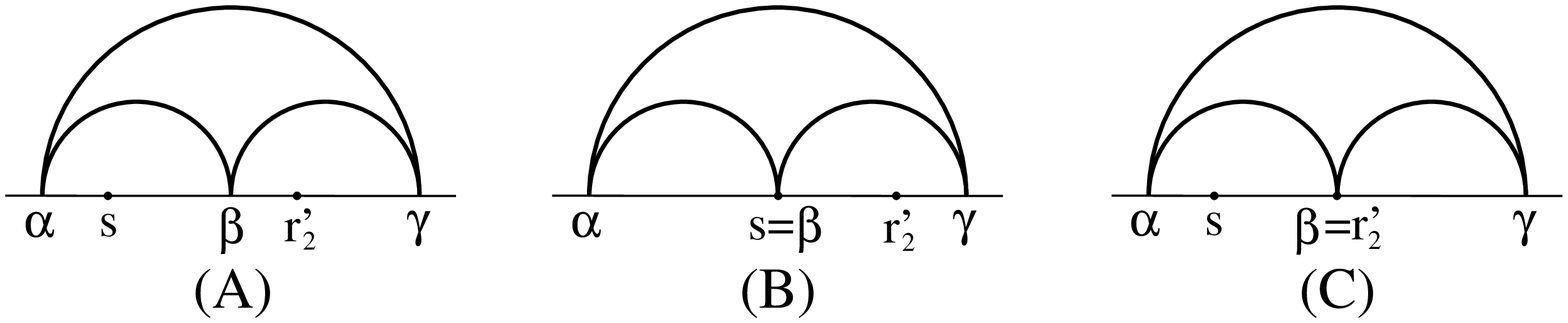}
\caption{The configurations of Farey arcs of Lemma~\ref{l:fareyarcs1}.}
\label{f:conf1}
\end{figure}

\begin{lem}\label{l:fareyarcs1}
Let $s, r'_2\in (-\infty,-1)$, with $s<r'_2$. Then, there exist 
a Farey arc \arc{$\al\ga$} with middle point $\be=m(\al,\ga)$ satisfying 
\[
-\infty\leq \al < s\leq \be\leq r'_2 < \ga\leq -1.
\]
In other words, there is a configuration of Farey arcs as in 
Figure~\ref{f:conf1}(A), \ref{f:conf1}(B) or \ref{f:conf1}(C).
\end{lem}

\begin{proof}
We provisionally define $\al=-\infty$ and  $\ga=-1$. If $s\leq \be=m(\al,\ga)\leq r'_2$ then 
$\al$, $\be$ and $\ga$ already satisfy the statement and the lemma is proved. 
Otherwise we have either $s>\be$ or $r'_2<\be$ (but not both, because $s<r'_2$). 
If $s>\be$ we redefine $\al=\be$, while if $r'_2<\be$ we redefine $\ga=\be$, 
and in both cases we set $\be$ equal to the new middle point $m(\al,\ga)$. 
As before, if $s\leq \be=m(\al,\ga)\leq r'_2$ we are done, otherwise 
either $s>\be$ or $r'_2<\be$ (but not both). Continuing in this fashion, 
after a finite number of steps we necessarily arrive at a 
configuration satisfying the statement of the lemma.
\end{proof}

\begin{figure}[!ht]
\centering
\includegraphics[scale=0.5]{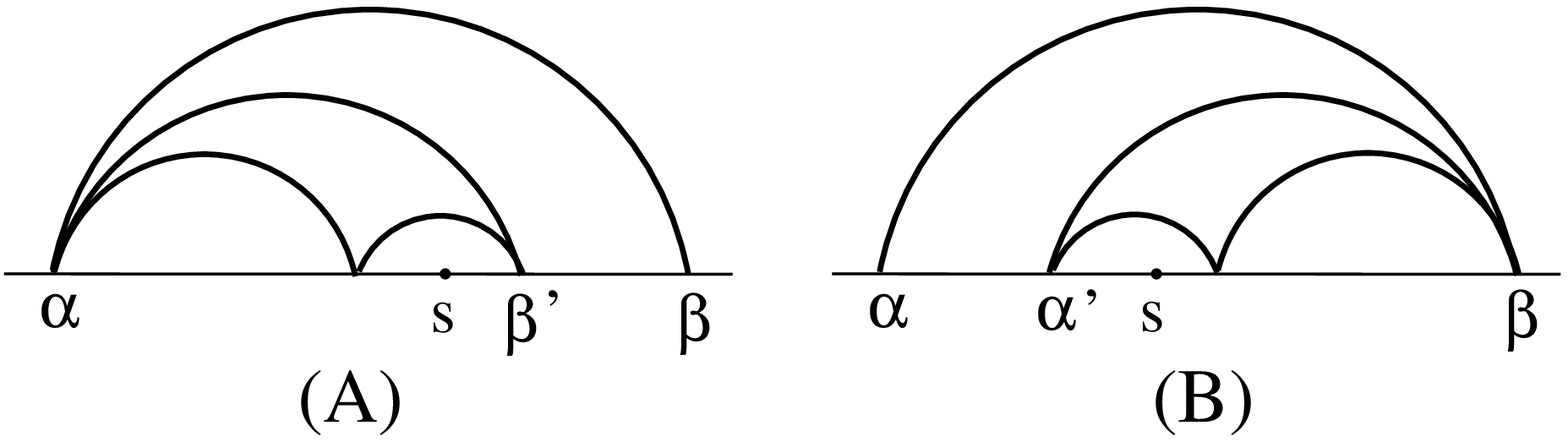}
\caption{The configurations of Farey arcs of Lemma~\ref{l:fareyarcs2}.}
\label{f:conf2}
\end{figure}

\begin{lem}\label{l:fareyarcs2}
Let \arc{$\al\be$} be a Farey arc with $\al,\be\in [-\infty,-1]$, and  
let $s\in (\al,\be)$. Then, there exist Farey arcs \arc{$\al\be'$} and 
\arc{$\al'\be$} such that 
\[
m(\al,\be')\leq s < \be'\leq \be\quad\text{and}\quad \al\leq \al' < s\leq m(\al',\be).
\]
In other words, there are configurations of Farey arcs as in 
Figure~\ref{f:conf2}(A) and~\ref{f:conf2}(B).
\end{lem}

\begin{proof} 
We only prove the existence of the arc \arc{$\al\be'$}. The existence of the arc \arc{$\al'\be$}
can be established in the same way. If $m(\al,\be)\leq s$ we define $\be'=\be$, and the 
statement is proved. If $s<m(\al,\be)$, we define (temporarily) 
$\be'=m(\al,\be)$. If $m(\al,\be')\leq s$ the arc \arc{$\al\be'$} satisfies the statement. 
Otherwise $s<m(\al,\be')$, we redefine $\be'=m(\al,\be')$ and we keep going in the same way. 
In a finite number of steps we are bound to find the arc \arc{$\al\be'$} with the 
stated property.  
\end{proof}

\begin{figure}[!ht]
\centering
\includegraphics[scale=0.5]{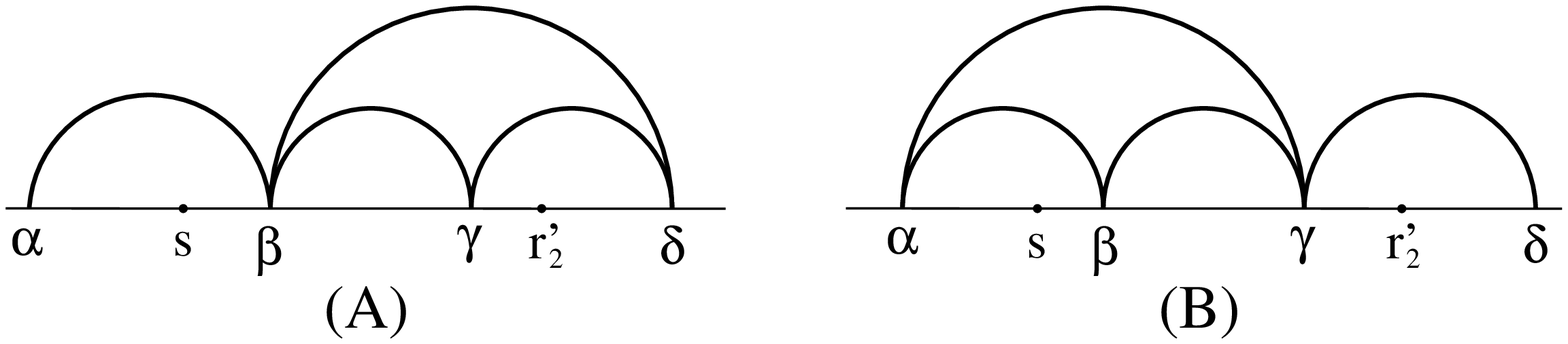}
\caption{The configurations of Farey arcs of Lemma~\ref{l:fareyarcs3}.}
\label{f:conf3}
\end{figure}

\begin{lem}\label{l:fareyarcs3} 
Let $s, r'_2\in (-\infty,-1)$, with $s<r'_2$. Then, there exist 
Farey arcs \arc{$\al\be$}, \arc{$\be\ga$}, \arc{$\ga\de$}, \arc{$xy$},
with $\text{\arc{$xy$}}\in\{\text{\arc{$\al\ga$}},\text{\arc{$\be\de$}}\}$, satisfying  
\[
-\infty\leq \al<s\leq \be<\ga\leq r'_2<\de\leq -1.
\]
In other words, there is a configuration of Farey arcs 
as in Figure~\ref{f:conf3}(A) or~\ref{f:conf3}(B).  
\end{lem} 

\begin{proof} 
By Lemma~\ref{l:fareyarcs1} there is a configuration of Farey arcs as in 
Figure~\ref{f:conf1}(A), \ref{f:conf1}(B) or \ref{f:conf1}(C).
If~\ref{f:conf1}(A) holds we can apply Lemma~\ref{l:fareyarcs2} to the point 
$s$ and the arc \arc{$\al\be$} of Figure~\ref{f:conf1}(A) to  
find a configuration as in Figure~\ref{f:conf2}(B). If we set $m=m(\al',\be)$, 
the Farey arcs \arc{$\al'm$}, \arc{$m\be$}, \arc{$\al'\be$} and \arc{$\be\ga$} 
provide a configuration as in Figure~\ref{f:conf3}(B). 
If~\ref{f:conf1}(B) holds we can apply Lemma~\ref{l:fareyarcs2} to the point 
$r'_2$ and the arc \arc{$\be\ga$} of Figure~\ref{f:conf1}(B) to find a configuration 
as in Figure~\ref{f:conf2}(A). In other words, there exists a 
Farey arc \arc{$\be\be'$} such that $m(\be,\be')\leq r'_2<\be'$. Setting $m=m(\be,\be')$, 
the Farey arcs \arc{$\al\be$}, \arc{$\be m$}, \arc{$m\be'$} and \arc{$\be\be'$} 
provide a configuration as in Figure~\ref{f:conf3}(A). 
Finally, if~\ref{f:conf1}(C) holds we can apply Lemma~\ref{l:fareyarcs2} to the point 
$s$ and the arc \arc{$\al\be$} of Figure~\ref{f:conf1}(C) to  
find a configuration as in Figure~\ref{f:conf2}(B). If we set $m=m(\al',\be)$, 
the Farey arcs \arc{$\al' m$}, \arc{$m\be$}, \arc{$\al'\be$} and \arc{$\be\ga$} 
provide a configuration as in Figure~\ref{f:conf3}(B). 
\end{proof}

\begin{thm}\label{t:r1r2big}
Suppose that $k\geq 3$, $1>r_1\geq r_2\geq \cdots \geq r_k>0$ and 
$r_1+r_2\geq 1$. Then, $Y(-1;r_1,\ldots,r_k)$ is orientation preserving 
diffeomorphic to the boundary of a Stein surface. 
\end{thm}

\begin{proof} 
It is easy to check 
that the condition $r_1+r_2\geq 1$ is equivalent to $s\leq r'_2$. If $s=r'_2$ Condition~\eqref{e:Gcond}
is automatically satisfied, therefore we may assume $s<r'_2$. By 
Lemma~\ref{l:fareyarcs3} there is a configuration of Farey arcs 
as in Figure~\ref{f:conf3}(A) or~\ref{f:conf3}(B). Let us suppose that the first 
case occurs. In view of the action of $PSL(2,\Z)$ on the Farey tessellation, there 
is a unique map $A\co\Q\cup\{\infty\}\to\Q\cup\{\infty\}$ of the form 
$A(r)=\frac{c+dr}{a+br}$ satisfying Condition~\eqref{e:g1} and such that 
$A(\be)=0$, $A(\ga)=\infty$ and $A(\de)=-1$. By construction, $A$ satisfies 
Conditions~\eqref{e:g2} and~\eqref{e:g3} as well, 
$c/a=A(0)\in (-1,0)$ and $A(s)\in (-1,0)$. According to Gompf's 
condition~\eqref{e:Gcond}, in order to prove 
that $Y(-1;r_1,\ldots,r_k)$ carries Stein fillable contact structures it suffices 
to show that 
\[
n_A(r'_1,r'_2) = -m(\llbracket t\rrbracket + 1) - M\geq -1, 
\]
where $t=1/A(s)\in (-\infty,-1)$, $M=\max(|a|,|c|)$ and $m=\min(|a|,|c|)$.
Observe that, since $A(0)=c/a\in (-1,0)$, $M=|a|$ and $m=|c|$. 
Condition~\eqref{e:g1} implies that if $m=0$ then 
$M=1$, therefore we may assume without loss of generality that $m > 0$. 
An easy calculation shows that the condition $-m(\llbracket t\rrbracket + 1) - M\geq -1$ 
is equivalent to 
\[
 \llbracket t\rrbracket + 1 - \frac1m \leq -\frac{M}{m} = \frac{1}{A(0)}.
\]
In order to prove this inequality it suffices to show that there is an integer $N$ 
strictly greater than $1/A(s)$ and less than or equal to $1/A(0)$, i.e.~such that 
\[
\llbracket t\rrbracket + 1 \leq N \leq \frac1{A(0)}.
\]
This condition is satisfied if and only if there exists a Farey arc~\arc{$\infty x$} 
with $1/A(s)<x\leq 1/A(0)$ or, equivalently, if and only if there exists a 
Farey arc~\arc{$y0$} with $A(0)\leq y<A(s)$. Setting $y:=A(\al)$, such an arc 
is provided by the image under $A$ of the Farey arc~\arc{$\al\be$} from 
Lemma~\ref{l:fareyarcs3}, because by construction $A(\be)=0$ and 
$0<\infty=-\infty\leq \al < s$. This concludes the proof under the assumption 
that when at the beginning of the argument we apply Lemma~\ref{l:fareyarcs3} 
we end up with a configuration of Farey arcs 
as in Figure~\ref{f:conf3}(A). In case the configuration is the one 
given by Figure~\ref{f:conf3}(B) we can argue in a similar way, so we just 
describe the steps where there is a difference. We choose the unique map $A$ 
such that $A(\al)=-1$, $A(\be)=0$ and $A(\ga)=+\infty=-\infty$. Then, 
$c/a=A(0)\in (-\infty,-1)$, $t=A(r'_2)\in (-\infty,-1)$, $M=|c|$, $m=|a|$ 
and as before we may assume without loss that $m > 0$. The same calculation 
as in the previous case shows that Gompf's condition is equivalent to 
\[
\llbracket A(r'_2)\rrbracket + 1 - \frac1m \leq -\frac{M}{m} = A(0).
\]
As before, this condition is satisfied if there exists a Farey 
arc -\arc{$\infty z$} with $A(r'_2)< z\leq A(0)$. Setting  $z:=A(\de)$, such 
an arc is provided by the image under $A$ of the 
Farey arc~\arc{$\ga\de$} from Lemma~\ref{l:fareyarcs3}, 
because by construction $A(\ga)=-\infty$ and $r'_2<\de\leq -1<0$.
\end{proof}

We are now ready to prove Theorem~\ref{t:mainexistence}. We restate the result 
for the reader's convenience:

\noindent{\bf Theorem~1.3}~
Let $Y$ be a closed, oriented, Seifert fibered 3--manifold which is not of special type. 
Then, $Y$ is orientation preserving diffeomorphic to the boundary of a Stein surface.  

\begin{proof}
Gompf showed~\cite[Theorem~5.4]{Go} that a closed, oriented, Seifert fibered 3--manifold $Y$
admits a Stein filling unless it is of the form $Y(-1;r_1,\ldots, r_k)$ 
with $k\geq 3$ and $1>r_1\geq r_2\geq\cdots\geq r_k>0$, and by~\cite{EHN, Na94} 
(see~\cite[pp.~157--158]{Na94}), if $r_1+\cdots + r_k \leq 1$ then 
$(r_1,\ldots, r_k)$ is realizable. Thus, applying~Theorem~\ref{t:realizable} we 
conclude that if $Y$ is not of special type then $Y$ is diffeomorphic 
to the boundary of a Stein surface unless $Y$ is of the form 
$Y(-1;r_1,\ldots, r_k)$ with $k\geq 3$, $1>r_1\geq r_2\geq\cdots\geq r_k>0$ and 
$r_1+r_2\geq 1$. But the latter conditions are precisely the assumptions of 
Theorem~\ref{t:r1r2big}, therefore if $Y$ is not of special type then $Y$ is 
necessarily the boundary of a Stein surface.
\end{proof} 

\section{Nonexistence of symplectic fillings}\label{s:nonexistence}

The purpose of this section is to establish Theorem~\ref{t:mainnonexistence}. We start with 
some preliminaries, then we prove three auxiliary lemmas. After that we prove 
Theorem~\ref{t:mainnonexistence}.

\begin{figure}[ht]
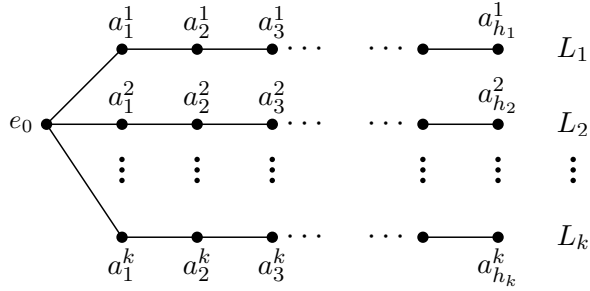

\begin{center}
\begin{graph}(10,3.5)(-9,-2)
\graphnodesize{0.15}
  \roundnode{n0}(-8,0)
  \roundnode{n11}(-7,1)
  \roundnode{n12}(-6,1)
  \roundnode{n13}(-5,1)
  \roundnode{n14}(-3,1)
  \roundnode{n15}(-2,1)
  
  \roundnode{n21}(-7,0)
  \roundnode{n22}(-6,0)
  \roundnode{n23}(-5,0)
  \roundnode{n24}(-3,0)
  \roundnode{n25}(-2,0)

  \roundnode{nh1}(-7,-1.5)
  \roundnode{nh2}(-6,-1.5)
  \roundnode{nh3}(-5,-1.5)
  \roundnode{nh4}(-3,-1.5)
  \roundnode{nh5}(-2,-1.5)

  \edge{n0}{n11}
  \edge{n11}{n12}
  \edge{n12}{n13}
  \edge{n14}{n15}
 
  \edge{n0}{n21}
  \edge{n21}{n22}
  \edge{n22}{n23}
  \edge{n24}{n25}

  \edge{n0}{nh1}
  \edge{nh1}{nh2}
  \edge{nh2}{nh3}
  \edge{nh4}{nh5}

  \autonodetext{n0}[w]{\small $e_0$}
  \autonodetext{n11}[n]{$a^1_1$}
  \autonodetext{n12}[n]{$a^1_2$}
  \autonodetext{n13}[n]{$a^1_3$}
  \autonodetext{n13}[e]{\large $\cdots$}
  \autonodetext{n14}[w]{\large $\cdots$}
  \autonodetext{n15}[n]{$a^1_{h_1}$}

  \autonodetext{n21}[n]{$a^2_1$}
  \autonodetext{n22}[n]{$a^2_2$}
  \autonodetext{n23}[n]{$a^2_3$}
  \autonodetext{n23}[e]{\large $\cdots$}
  \autonodetext{n24}[w]{\large $\cdots$}
  \autonodetext{n25}[n]{$a^2_{h_2}$}

  \autonodetext{nh1}[s]{$a^k_1$}
  \autonodetext{nh2}[s]{$a^k_2$}
  \autonodetext{nh3}[s]{$a^k_3$}
  \autonodetext{nh3}[e]{\large $\cdots$}
  \autonodetext{nh4}[w]{\large $\cdots$}
  \autonodetext{nh5}[s]{$a^k_{h_k}$}

  \freetext(-1,1){$L_1$}
  \freetext(-1,0){$L_2$}
  \freetext(-1,-1.5){$L_k$}

  \freetext(-7,-0.5){\Huge $\vdots$}
  \freetext(-6,-0.5){\Huge $\vdots$}
  \freetext(-5,-0.5){\Huge $\vdots$}
  \freetext(-3,-0.5){\Huge $\vdots$}
  \freetext(-2,-0.5){\Huge $\vdots$}
  \freetext(-1,-0,5){\Huge $\vdots$}
\end{graph}
\end{center}
\caption{The weighted star--shaped graph $\Ga$}
\label{f:graph} 
\end{figure}

Let $Y=Y(e_0;r_1,...,r_k)$ denote the oriented, Seifert fibered $3$--manifold given 
by the surgery description of Figure~\ref{f:seifert}, where 
$e_0\in\Z$, $r_i\in(0,1)\cap\Q$ and $r_1\geq r_2\geq\cdots\geq r_k$. 
The oriented 3--manifold $Y$ is the oriented boundary of the 4--dimensional 
plumbing $P_\Ga$ of $D^2$--bundles over 2--spheres described by the star--shaped 
weighted graph $\Gamma$ with $k$ legs illustrated in Figure~\ref{f:graph}. 

The weights of the vertices in the leg $L_i$, which form the 
string $(a_1^i,...,a_{h_i}^i)$, are given by the unique continued 
fraction expansion 
\[
-\frac{1}{r_i}=[a_1^i,...,a_{h_i}^i]:=a_1^i-\frac{1}{\displaystyle a_2^i - 
\frac{ \bigl. 1}{\displaystyle \ddots\ _{\displaystyle {a_{h_i-1}^i -\frac
{\bigl. 1}{a_{h_i}^i}}}} },\quad 
i=1,\ldots,k
\]
such that $a_j^i\leq -2$ for every $j$. 

We can associate to $\Gamma$ the intersection lattice $(\Z^{|\Gamma|},Q_\Gamma)$ of the plumbing 
$P_\Ga$. In the proof of Theorem~\ref{t:mainnonexistence} 
we will show that if $Y$ admits a symplectic filling then the intersection lattice of 
the plumbing associated to $-Y$ admits an isometric embedding into a 
standard diagonal lattice. Our present aim will be to prepare the ground 
for Lemma~\ref{l:no_e}, which shows that under a certain assumption such an 
embedding does not exist. 

We shall need Riemenschneider's point rule~\cite{Ri}, which we now recall. 
Let $p>q>0$ be coprime integers, and suppose 
\[
-\frac pq = [a_1,\ldots,a_l],\ a_i\leq -2,\quad 
-\dfrac p{p-q} = [b_1,\ldots,b_m],\ b_j\leq -2.
\]
Then, the coefficients $a_1,\ldots,a_l$ and $b_1,\ldots,b_m$ are related by 
a diagram of the form
\begin{align*}
\bullet \cdots\cdots \bullet\hspace{6cm}\\
\hspace{1.3cm}\bullet \cdots\cdots\bullet \hspace{4.7cm}\\
\hspace{1.8cm}\ddots\hspace{4.2cm} \\
\hspace{3.4cm}\bullet\cdots\cdots \bullet\hspace{2.6cm} \\
\hspace{4.7cm}\bullet \cdots\cdots \bullet\hspace{1.3cm}
\end{align*}
where the $i$--th row contains $|a_i|-1$ ``points" for $i=1,\ldots,l$, and the first point of each row is 
vertically aligned with the last point of the previous row. The point rule says that there are
$m$ columns, and the $j$--th column contains $|b_j|-1$ points for every $j=1,\ldots, m$.
For example if $-7/5=[-2,-2,-3]$ and $-7/2 =[-4,-2]$ the corresponding 
diagram is given by 
\begin{align*}
\bullet & \\
\bullet & \\
\bullet & \quad \bullet
\end{align*}

Let $\Ga$ be either a star--shaped or a linear weighted graph. If $(\Z^{|\Gamma|},Q_\Gamma)$ 
admits an embedding 
into a standard diagonal lattice $(\Z^k,-\mathrm{Id})$ with basis $e_1,...,e_k$, 
then we will write, for every subset $S$ of the set of vertices of $\Gamma$,
\[
U_S:=\{e_i\ |\ e_i\cdot v\neq 0\ \mbox{for some}\ v\in S\}.
\]

We can now start to work towards Lemma~\ref{l:no_e} and Theorem~\ref{t:mainnonexistence}. 

\begin{lem}\label{l:cuadrado}
Suppose $-\frac{1}{r}=[a_1,...,a_n]$ and  $-\frac{1}{s}=[b_1,...,b_m]$ where
$r+s=1$. Consider a weighted linear graph $\Psi$ 
having two connected components, $\Psi_1$ and $\Psi_2$, where $\Psi_1$ consists of
$n$ vertices $v_1,...,v_n$ with weights $a_1,...,a_n$ and $\Psi_2$ of $m$ vertices $w_1,...,w_m$
with weights $b_1,...,b_m$. 
Moreover, suppose that there is an embedding of the lattice $(\Z^{n+m},Q_\Psi)$ 
into $(\Z^k,-\mathrm{Id})$, with basis $e_1,...,e_k$, such that $e_1\in
U_{v_1}\cap U_{w_1}$ and $U_\Psi=\{e_1,...,e_k\}$. 
Then, 
\begin{itemize}
\item[$(1)$] $U_{\Psi_1}=U_{\Psi_2}$, and 
\item[$(2)$]$k=n+m$.
\end{itemize}
\end{lem}

\begin{proof}
We start showing that $(1)$ implies $(2)$. In fact, since
$r+s=1$, by \cite[Lemma~2.6]{Li07}, we have
$$\sum_{i=1}^n(-a_i-3)+\sum_{j=1}^m(-b_j-3)=-2$$
and therefore 
\begin{equation}\label{e:I}
|\mathrm{tr}\,(Q_\Psi)|=\sum |a_i|+\sum |b_j|=3(n+m)-2.
\end{equation}

The matrix $Q_\Psi$ is not singular \cite[Remark~2.1]{Li07} so we necessarily
have $k\geq n+m$. Let us write $k=n+m+x$ 
for some $x\geq 0$. Since we are assuming $(1)$, each vector of the basis
$e_1,...,e_k$ satisfies $e_i\in U_{\Psi_1}\cap U_{\Psi_2}$ 
and therefore $|\mathrm{tr}\,(Q_\Psi)|\geq 2k$. Moreover, since the graph $\Psi$
has $n+m-2$ edges, it follows that 
$|\mathrm{tr}\,(Q_\Psi)|\geq 2k+n+m-2=3(n+m)-2+2x$. Hence, by \eqref{e:I}, $x=0$
and $(2)$ holds.

We now assume that $(1)$ does not hold and show that this assumption leads to a
contradiction. Start by defining the 
sets $E_1:=U_{\Psi_1}\setminus U_{\Psi_2}$, $E_2:=U_{\Psi_2}\setminus
U_{\Psi_1}$ and $E_{12}:=U_{\Psi_1}\cap U_{\Psi_2}$. 
Since we are assuming that $(1)$ does not hold, we have $E_1\neq\emptyset$ or
$E_2\neq\emptyset$. By simmetry, we may assume that $E_1\neq\emptyset$. 
It follows that there exists a smallest
index $n_0\in\{1,...,n\}$ such that 
$U_{v_{n_0}}\cap E_1\neq\emptyset$. This condition allows us to construct a new
connected linear graph $\bar\Psi_1$ with $n_0$ 
vertices and associated string of weights $(a_1,...,a_{n_0-1},\bar a_{n_0})$, where 
$$\bar a_{n_0}:=a_{n_0}+\sum_{e_\ell\in E_1}(v_{n_0}\cdot e_\ell)^2.$$
Notice that since the intersection lattice associated with $\Psi$ admits an 
embedding into a diagonal lattice, there is a naturally induced analogous embedding 
of the intersection lattice associated with $\bar\Psi_1\cup\Psi_2$. We claim that 
$\bar a_{n_0}\leq -2$. In fact, if $n_0=1$ the assumption $e_1\in U_{v_1}\cap
U_{w_1}$ and the equality $v_1\cdot w_1=0$ 
imply $|E_{12}\cap U_{v_1}|\geq 2$ and therefore in this case $\bar a_1\leq -2$.
On the other hand, if $n_0>1$ then, by 
definition of $n_0$, it holds that $U_{v_{n_0-1}}\subseteq E_{12}$. The
equalities $v_{n_0-1}\cdot v_{n_0}=1$ and 
$v_{n_0}\cdot w_\ell=0$ for every $\ell\in\{1,...,m\}$ force $|E_{12}\cap
U_{v_{n_0}}|\geq 2$ and therefore 
$\bar a_{n_0}\leq -2$ as claimed.

Now, since $-2\geq\bar a_{n_0}> a_{n_0}$ we have,  by standard facts on
continued fractions, that 
$-\frac{1}{\bar r}:=[a_1,...,a_{n_0-1},\bar a_{n_0}]$ satisfies $\bar r+s>1$
(since $r+s=1$). 
Let $\bar r'$ be such that $\bar r+\bar r'=1$. Using Riemenschneider's point
diagram it is not difficult to 
check  that $r+s=1$ and $-2\geq\bar a_{n_0}>a_{n_0}$ imply that there is some
$t<m$  such that  
$-\frac{1}{\bar r'}=[b_1,...,b_t]$. Let us call $\bar\Psi'_1\subseteq\Psi_2$ the
linear connected subgraph with 
associated string of weights $(b_1,...,b_t)$. There are two possibilities, either
$U_{\bar\Psi_1}=U_{\bar\Psi'_1}$ or 
$U_{\bar\Psi_1}\neq U_{\bar\Psi'_1}$. 

Notice that by construction $|\bar\Psi_1|\leq|\Psi_1|$ and
$|\bar\Psi'_1|<|\Psi_2|$. Moreover, if 
$|\bar\Psi_1|=1$ [resp.\ $|\bar\Psi'_1|=1$] then $\bar\Psi'_1$ [resp.\
$\bar\Psi_1$] is a $(-2)$--chain and it is 
immediate to check that in this case, since $e_1\in U_{v_1}\cap U_{w_1}$, it
holds $U_{\bar\Psi_1}=U_{\bar\Psi'_1}$. 
Since $\bar r+\bar r'=1$, if $U_{\bar\Psi_1}\neq U_{\bar\Psi'_1}$ we can repeat
the above construction with 
$\bar\Psi_1$ and $\bar\Psi'_1$ playing the role of $\Psi_1$ and $\Psi_2$. It
follows that, after a finite number 
of steps, we necessarily obtain from $\Psi_1$ and $\Psi_2$ two linear weighted graphs,
which we still call $\bar\Psi_1$ and 
$\bar\Psi'_1$, such that $U_{\bar\Psi_1}=U_{\bar\Psi'_1}$ and either
$\bar\Psi_1\subseteq\Psi_1$ or 
$\bar\Psi'_1\subseteq\Psi_2$. By simmetry we may assume that
$\bar\Psi'_1\subseteq\Psi_2$. 

Since the strings of weights associated to $\bar\Psi_1$ and $\bar\Psi'_1$ are related to
one another by Riemenschneider's point 
rule, we know, by the first part of this proof and using the same notation, 
that $|U_{\bar\Psi_1}\cup
U_{\bar\Psi'_1}|=n_0+t$. Consider 
the vector $w_{t+1}$. Since $U_{\bar\Psi_1}=U_{\bar\Psi'_1}$, $w_t\cdot
w_{t+1}=1$ and $w_{t+1}\cdot v_\ell=0$ 
for every $\ell\in\{1,...,n_0\}$, the vector
$$\bar w_{t+1}:=w_{t+1}+\sum_{e_i\not\in U_{\bar\Psi_1}}(e_i\cdot
w_{t+1})\,e_i$$
satisfies $\bar w_{t+1}\cdot\bar w_{t+1}\leq -2$. It follows that the disconnected linear 
graph 
$\bar\Psi_1\cup\bar\Psi'_1\cup\{w_{t+1}\}$, which has $n_0+t+1$ vertices admits
an embedding into a diagonal 
lattice of rank $|U_{\Psi_1}|=n_0+t$ which contradicts
\cite[Remark~2.1]{Li07}.
\end{proof}

\begin{lem}\label{l:trunc}
Let $-\frac{1}{r}=[a_1,...,a_n]$ and $-\frac{1}{s}=[b_1,...,b_m]$ be such  that
$r+s>1$. Then there exists 
$n_0\leq n$ and $m_0\leq m$ such that $-\frac{1}{r_0}=[a_1,...,a_{n_0}]$ and
$-\frac{1}{s_0}=[b_1,...,b_{m_0}]$ 
satisfy $r_0+s_0=1$. 
\end{lem}

\begin{proof}
Let $r'$ be such that $r + r' =1$ and suppose 
$-1/r' =[a'_1,...,a'_{n'}]$. Since
$s> r'$, by standard facts on 
continued fractions there are two possibilities: either $b_i=a'_i$ for all
$i\in\{1,...,n'\}$ and $m>n'$ or 
there is a smallest index $k$ such that $b_k>a'_k$. In the first case we set
$n_0=n$ and $m_0=n'$. In the 
second case let us consider the first $k$ columns of dots in the
Riemenschneider's point diagram obtained 
from $(a_1,...,a_n)$. Then, $n_0$ equals the number of rows in this diagram
minus $b_k-a'_k$ and $m_0=k$. 
Note that in this way $[a_1,...,a_{n_0}]$ and $[b_1,...,b_{m_0}]$ are related to
one another by Riemenschneider's 
point rule and therefore $r_0+s_0=1$.
\end{proof}

\begin{lem}\label{l:no_e}
Suppose $k\geq 3$ and $1>r_1\geq\cdots\geq r_k>0$ and $r_{k-1}+r_k>1$. 
Then, the intersection lattice of the plumbing associated to $Y:=Y(-k+1;r_1,...,r_k)$ 
cannot be embedded into a negative diagonal standard lattice.
\end{lem}

\begin{proof}
Let $\Gamma$ be the plumbing graph of Figure~\ref{f:graph} associated to $Y$, 
and suppose by contradiction that there exists an embedding 
of $(\Z^{|\Gamma|},Q_\Gamma)$ into $(\Z^d,-\mathrm{Id})$ with basis
$e_1,...,e_d$ for some $d\geq |\Gamma|$. We will use the following notations 
for the vertices of $\Gamma$: $v_0$ for the central vertex and $v_j^i$ for the vertices in the legs, 
where $i$ indicates the leg to which $v_j^i$belongs and $j$ the position in the leg, 
with $j=1$ being the index of the vertex connected to the central vertex.

Since $\Gamma$ has $k$ legs connected to the central vertex which has weight
$-k+1$, there must exist some basis 
vector, say $e_1$, and two legs, say $L_i$ and $L_j$, such that the products
$v_0\cdot e_1$, $v_1^i\cdot e_1$ and 
$v_1^j\cdot e_1$ are not $0$.

Let $-\frac{1}{r_i}=[a_1^i,...,a_p^i]$ and $-\frac{1}{r_j}=[a_1^j,...,a_q^j]$.
Since $r_i+r_j>1$, by
Lemma~\ref{l:trunc} there exist $p_0\leq p$ and $q_0\leq q$ such that the
strings  $(a_1^i,...,a_{p_0}^i)$ and 
$(a_1^j,...,a_{q_0}^j)$ are related to one another by Riemenschneider's point
rule. Moreover, since 
$e_1\in U_{v_1^i}\cap U_{v_1^j}$ Lemma~\ref{l:cuadrado} applies and therefore
the disconnected subgraph 
$\Psi\subseteq\Gamma$ consisting of the vertices
$v_1^i,...,v_{p_0}^i,v_1^j,...,v_{q_0}^j$ satisfies $|U_{\Psi}|=p_0+q_0$. 
Furthermore, writing $\Psi=\Psi_1\cup\Psi_2$ where $\Psi_1$ [resp.\ $\Psi_2$]
consists of the vertices $v_1^i,...,v_{p_0}^i$ 
[resp.\ $v_1^j,...,v_{q_0}^j$] we have, by Lemma~\ref{l:cuadrado}\,$(1)$, 
$U_{\Psi_1}=U_{\Psi_2}$.

Now, since $r_i+r_j>1$ we have $(p_0,q_0)\neq (p,q)$ so we can assume without loss
of generality $p>p_0$. The 
vector $v_{p_0+1}^i$ satisfies $v_{p_0+1}^i\cdot v_{p_0}^i=1$ and the equality
$U_{\Psi_1}=U_{\Psi_2}$ implies 
that $|U_{\Psi}\cap U_{v_{p_0+1}}|\geq 2$. Consider the vector
$$\bar v_{p_0+1}:=v_{p_0+1}+\sum_{e_i\not\in U_{\Psi}}(e_i\cdot
v_{p_0+1})\,e_i,$$
which by construction satisfies $U_{\bar v_{p_0+1}}\subseteq U_\Psi$ and
$v_{p_0+1}\cdot v_{p_0+1}\leq -2$. It 
follows that the linear graph $\Psi\cup\{v_{p_0+1}\}$, which has $p_0+q_0+1$
vertices admits an embedding into a 
diagonal lattice of rank $|U_{\Psi}|=p_0+q_0$ which contradicts
\cite[Remark~2.1]{Li07}.
\end{proof}

We are now ready to prove Theorem~\ref{t:mainnonexistence}. We restate the 
result for the reader's convenience:

\noindent{\bf Theorem~1.4}~
A closed, oriented, Seifert fibered 3--manifold of special type admits 
no symplectic fillings. 

\begin{proof}
Suppose that the oriented 3--manifold $Y$ is orientation 
preserving diffeomorphic to $Y(-1;r_1,\ldots, r_k)$, where $k\geq 3$, 
$1>r_1\geq r_2\geq \cdots \geq r_k>0$ and $(r_1,\ldots,r_k)$ satisfies Conditions (1) 
and (2) of Definition~\ref{d:specialtype}. The fact that $(r_1,\ldots,r_k)$ is not realizable 
implies that $Y$ is an $L$--space~\cite{LS07-2}. Therefore, if $Y$ admits a symplectic filling $W$ then 
$b_2^+(W)=0$~\cite[Theorem~1.4]{OS04}. 
Consider the space $-Y=Y(1-k;\overline{r_1},...,\overline{r_k})$ where 
$\overline{r_i}:=1-r_{k-i+1}$ for $i\in\{1,...,k\}$. Since 
\[
e(-Y) := 1-k + \sum_i \overline r_i = 1 - \sum_i r_i <0, 
\]
by~\cite[Theorem~5.2]{NR78} there is a negative definite plumbing graph $\Gamma$ such that 
$-Y=\partial P_\Gamma$. Consider the $4$--manifold $X$ obtained gluing together 
$P_\Gamma$ and $W$ along their common boundary. By construction $X$ is smooth, 
closed and negative, therefore by Donaldson's celebrated theorem~\cite{Do87} 
its associated intersection form must be diagonalizable. It follows that 
if $Y$ admitted a symplectic filling then the intersection lattice of $P_\Ga$ 
would admit an embedding into a diagonal, negative standard lattice. 
The assumption $1>r_1+r_2$ for $Y$ reads
$\overline{r_{k-1}}+\overline{r_k}>1$ for $-Y$, which by Lemma~\ref{l:no_e} 
implies that the intersection lattice of $P_\Ga$ does not admit an embedding into a
diagonal lattice. Therefore we conclude that $Y$ admits no 
symplectic fillings. 
\end{proof}


\bibliographystyle{plain}
\bibliography{biblio}

\begin{thebibliography}{10}

\bibitem{EHN}
U.~Hirsch {D. Eisenbud} and W.~Neumann.
\newblock Transverse foliations of {S}eifert bundles and self-homeomorphism of
  the circle.
\newblock {\em Comment. Math. Helv.}, 56(4):638--660, 1981.

\bibitem{Do87}
S.~Donaldson.
\newblock The orientation of {Y}ang-{M}ills moduli spaces and 4-manifold
  topology.
\newblock {\em J. Differ. Geom.}, 26:397--428, 1987.

\bibitem{El1}
Y.~Eliashberg.
\newblock Topological characterization of {S}tein manifolds of dimension $> 2$.
\newblock {\em Internat. J. of Math.}, 1:29--46, 1990.

\bibitem{El2}
Y.~Eliashberg.
\newblock Fillings by holomorphic discs and its applications.
\newblock {\em London Math. Soc. Lecture Notes Series}, 151:45--67, 1991.

\bibitem{Et}
J.~B. Etnyre.
\newblock Introductory lectures on contact geometry.
\newblock In {\em Topology and geometry of manifolds ({A}thens, {GA}, 2001)},
  volume~71 of {\em Proc. Sympos. Pure Math.}, pages 81--107, Providence, RI,
  2003. Amer. Math. Soc.

\bibitem{Ge}
H.~Geiges.
\newblock {\em An introduction to contact topology}, volume 109 of {\em
  Cambridge Studies in Advanced Mathematics}.
\newblock Cambridge University Press, Cambridge, 2008.

\bibitem{Gh}
P.~Ghiggini.
\newblock Strongly fillable contact 3-manifolds without {S}tein fillings.
\newblock {\em Geom. Topol.}, 9:1677--1687, 2005.

\bibitem{GLS}
P.~Ghiggini, P.~Lisca, and A.~I. Stipsicz.
\newblock Tight contact structures on some small {S}eifert fibered 3-manifolds.
\newblock {\em Amer. J. Math.}, 129(5):1403--1447, 2007.

\bibitem{Go}
R.~Gompf.
\newblock Handlebody construction of {S}tein surfaces.
\newblock {\em Ann. Math.}, 148:619--693, 1998.

\bibitem{Gr}
M.~Gromov.
\newblock Pseudo--holomorphic curves in symplectic manifolds.
\newblock {\em Invent. Math.}, 82:307--347, 1985.

\bibitem{Li98}
P.~Lisca.
\newblock {Symplectic fillings and positive scalar curvature.}
\newblock {\em Geom. Topol.}, 2:103--116, 1998.

\bibitem{Li07}
P.~Lisca.
\newblock {Lens spaces, rational balls and the ribbon conjecture.}
\newblock {\em Geom. Topol.}, 11:429--472, 2007.

\bibitem{LS03}
P.~Lisca and A.~I. Stipsicz.
\newblock An infinite family of tight, not semi-fillable contact
  three-manifolds.
\newblock {\em Geom. Topol.}, 7:1055--1073, 2003.

\bibitem{LS04-2}
P.~Lisca and A.~I. Stipsicz.
\newblock Ozsv{\'a}th-{S}zab{\'o} invariants and tight contact three-manifolds.
  {I}.
\newblock {\em Geom. Topol.}, 8:925--945, 2004.

\bibitem{LS04-1}
P.~Lisca and A.~I. Stipsicz.
\newblock Tight, not semi-fillable contact circle bundles.
\newblock {\em Math. Ann.}, 328(1-2):285--298, 2004.

\bibitem{LS07-1}
P.~Lisca and A.~I. Stipsicz.
\newblock Ozsv{\'a}th-{S}zab{\'o} invariants and tight contact three-manifolds.
  {II}.
\newblock {\em J. Differential Geom.}, 75(1):109--141, 2007.

\bibitem{LS09}
P.~Lisca and A.~I. Stipsicz.
\newblock On the existence of tight contact structures on {S}eifert fibered
  3-manifolds.
\newblock {\em Duke Math. J.}, 148(2):175--209, 2009.

\bibitem{LS07-2}
P.~Lisca and A.I. Stipsicz.
\newblock Ozsv{\'a}th-{S}zab{\'o} invariants and tight contact 3-manifolds.
  {III}.
\newblock {\em J. Symplectic Geom.}, 5(4):357--384, 2007.

\bibitem{LP01}
A.~Loi and R.~Piergallini.
\newblock Compact {S}tein surfaces with boundary as branched covers of {$B^4$}.
\newblock {\em Invent. Math.}, 143(2):325--348, 2001.

\bibitem{Na94}
R.~Naimi.
\newblock Foliations transverse to fibers of {S}eifert manifolds.
\newblock {\em Comment. Math. Helv.}, 69(1):155--162, 1994.

\bibitem{NR78}
W.~D. Neumann and F.~Raymond.
\newblock {S}eifert manifolds, plumbing, $\mu$-invariant and orientation
  reversing maps.
\newblock Algebraic and geometric topology (Santa Barbara 1977), {\it Lecture
  Notes in Math.}, 664:163--196 (1978).

\bibitem{OS04}
P.~Ozsv{\'a}th and Z.~Szab{\'o}.
\newblock {Holomorphic disks and genus bounds.}
\newblock {\em Geom. Topol.}, 8:311--334, 2004.

\bibitem{OS05}
P.~Ozsv{\'a}th and Z.~Szab{\'o}.
\newblock On knot {F}loer homology and lens space surgeries.
\newblock {\em Topology}, 44(6):1281--1300, 2005.

\bibitem{Ri}
O.~Riemenschneider.
\newblock Deformationen von {Q}uotientensingularit{\"a}ten (nach zyklischen
  {G}ruppen).
\newblock {\em Math. Ann.}, 209:211--248, 1974.

\end{thebibliography}

\end{document}